\documentclass[12pt]{amsart}
\usepackage[all]{xy}
\usepackage{times}
\usepackage{amscd,amssymb,graphics}
\usepackage[T1]{fontenc}
\usepackage[english,french]{babel}
\usepackage[ansinew]{inputenc}
\usepackage[all]{xy}
\usepackage[centertags]{amsmath}
\usepackage{latexsym}
\usepackage{amsfonts}
\usepackage{amssymb}
\newtheorem{theorem}{Th\'eor\`eme}[section]
\newtheorem{corollary}[theorem]{Corollaire}
\newtheorem{lemma}[theorem]{Lemme}
\newtheorem{proposition}[theorem]{Proposition}
\theoremstyle{definition}
\newtheorem{definition}[theorem]{D\'efinition}
\newtheorem{question}[theorem]{Question}

\newtheorem{remark}[theorem]{Remarque}

\numberwithin{equation}{section}

\begin{document}

\title[Sur les espaces test pour la moyennabilit\'e]{Sur les espaces test pour la moyennabilit\'e}

\author[Y. Al-Gadid]{Yousef Al-Gadid}
\address{Department of Mathematics, Faculty of Science, Al-Fateh University, Tripoli-Lybie.}
\email{yousef$\_$algadid@yahoo.com}

\author[B. R. MBOMBO]{Brice R. MBOMBO}
\address{D\'epartement de Math\'ematiques, Facult\'e des Sciences,
Universit\'e de Yaound\'e I, BP 812 Yaound\'e, Cameroun.}
\email{bricero@yahoo.fr}

\author[V. G. Pestov]{Vladimir G. Pestov}
\address{Département de Math\'ematiques et Statistiques,
Universit\'e  d'Ottawa, 585, ave. King Edward, Ottawa, Ontario, Canada K1N 6N5.}
\email{vpest283@uottawa.ca}


\begin{abstract}
Nous observons qu'un groupe polonais $G$ est moyennable si et seulement si toute action continue de $G$ sur le cube de Hilbert poss\`ede une mesure de probabilit\'e invariante. Cela g\'en\'eralise un r\'esultat de Bogatyi et Fedorchuk. Nous d\'emontrons \'egalement que les actions continues sur l'espace de Cantor permettent de tester la moyennabilit\'e, la moyennabilit\'e extr\^eme des groupes polonais non archimediens, et la moyennabilit\'e \`a l'infini des groupes discrets d\'enombrables. Il en r\'esulte que cette derni\`ere propri\'et\'e peut \'egalement \^etre test\'ee par les actions sur le cube de Hilbert. Ces r\'esultats g\'en\'eralisent un crit\`ere de Giordano et de la Harpe.
\\[3mm]
{\it Abstract:}
We observe that a Polish group $G$ is amenable if and only if every continuous action of $G$ on the Hilbert cube admits an invariant probability measure. This generalizes a result of Bogatyi and Fedorchuk. We also show that actions on the Cantor space can be used to detect amenability and extreme amenability of Polish nonarchimedean groups as well as amenability at infinity of discrete countable groups. As corollary, the latter property can also be tested by actions on the Hilbert cube. These results generalise a criterion due to Giordano and de la Harpe.
\end{abstract}

\maketitle

\section{Introduction}
Un groupe topologique $G$ est dit {\em moyennable} si toute action continue de $G$ sur un espace compact $X$ poss\`ede une mesure de probabilit\'e bor\'elienne invariante. En r\'eponse \`a une question de Grigorchuk, Giordano et de la Harpe \cite{gio} ont montr\'e  qu'un groupe discret d\'enombrable $G$ est moyennable si et seulement si toute action continue de $G$ sur l'ensemble de Cantor $D^{\aleph_{0}}$ poss\`ede une mesure de probabilit\'e invariante. On peut dire que l'ensemble de Cantor est un {\em espace test} pour la moyennabilit\'e des groupes discrets d\'enombrables.
Dans le m\^eme sens, Bogatyi et Ferdorchuk \cite{boga} ont r\'epondu \`a une question de \cite{gio} et demontr\'e que le cube de Hilbert $I^{\aleph_{0}}$ est \'egalement un {\em espace test} pour la moyennabilit\'e des groupes discrets d\'enombrables.

Dans cet article, nous d\'emontrons que le cube de Hilbert reste un espace test pour la moyennabilit\'e de tous les groupes polonais. Nous d\'emontrons \'egalement que le r\'esultat de Giordano et de la Harpe reste vrai pour les groupes polonais non archim\'ediens.

Un groupe topologique $G$ est dit {\em extr\^emement moyennable} si toute action continue de $G$ sur un espace compact poss\`ede un point fixe. Nous d\'emontrons que l'ensemble de Cantor est un espace test pour la moyennabilit\'e extr\^eme des groupes polonais non archim\'ediens.
La question d'existence d'un espace test pour la moyennabilit\'e extr\^eme des groupes polonais reste ouverte, car le cube de Hilbert ne poss\`ede pas cette propri\'et\'e.

Si un groupe discret d\'enombrable $G$ op\`ere par hom\'eomorphismes sur un espace compact $X$, alors l'action de $G$ sur $X$ est {\em moyennable} s'il existe une suite $b^{n}:X\longrightarrow \mathcal{P}(G)$ d'applications continues tel que: $\underset{n\longrightarrow \infty}{\overset{}{\lim}}\,\,\underset{x\in X}{\overset{}{\sup}}\|gb^{n}_{x}-b^{n}_{gx}\|_{1}=0$ pour tout $g\in G$, o\`u $\mathcal{P}(G)$ d\'esigne l'espace des mesures de probabilit\'e sur $G$, muni de la topologie vague.
Un groupe discret d\'enombrable $G$ est dit {\em moyennable \`a l'infini} \cite{higson},\cite{anan}, ou {\em topologiquement moyennable}, s'il existe un espace compact $X$ et une action par hom\'eomorphismes de $G$ sur $X$ qui est moyennable. Les exemples des groupes moyennables \`a l'infini comprennent les groupes moyennables, les groupes classiques sur les corps, les groupes d'automorphismes d'arbres enracin\'es r\'eguliers, les groupes hyperboliques (en particulier, les groupes libres). Voir \cite{anan} et \cite{BO}, Ch. 5.

Par analogie avec le r\'esultat de Giordano et de la Harpe, nous d\'emontrons qu'un groupe discret d\'enombrable $G$ est moyennable \`a l'infini si et seulement si $G$ poss\`ede une action moyennable sur l'ensemble de Cantor ou sur le cube de Hilbert. Autrement dit, l'ensemble de Cantor et le cube de Hilbert sont des espaces test pour la moyennabilit\'e topologique des groupes discrets d\'enombrables.

\section{D\'ecomposition du compactifi\'e de Samuel $\mathcal{S}(G)$ en limite inverse}
Dans cette section, nous rassemblons quelques r\'esultats bien connus des sp\'ecialistes, mais donc les r\'ef\'erences sont souvent difficiles \`a retrouver.

Soit $X$ un espace compact. On note $C(X)$ l'ensemble des fonctions continues sur $X$ \`a valeurs r\'eelles et $C(X)^{\prime}$ le dual topologique de $C(X)$. L'espace $\mathcal{P}(X)$ des mesures de probabilit\'e sur $X$ est un sous-espace compact de $C(X)^{\prime}$ muni de la topologie vague. Si $X$ est un $G$-espace compact, alors l'action continue de $G$ sur $X$ peut \^etre prolong\'ee \`a $\mathcal{P}(X)$ de mani\`ere \'evidente.

Soient $G$ un groupe topologique et $E=RUCB(G)$ la $C^{\star}$-alg\`ebre commutative et unif\`ere constitu\'ee des fonctions $x:G\longrightarrow \mathbb{C}$ born\'ees et uniform\'ement continues \`a droite: pour tout $\varepsilon > 0$, il existe un voisinage $V$ de $e$ tel que $gh^{-1}\in V\Rightarrow |x(g)-x(h)|< \varepsilon$.

L'espace de Gelfand $\mathcal{S}(G)$ de $E$, c'est-\`a-dire l'espace de id\'eaux maximaux de $E$ muni de sa topologie compacte usuelle, est un compactifi\'e de $G$ \cite{brook}. C'est le compactifi\'e de Samuel \cite{samuel} de $G$ par rapport \`a la structure uniforme droite sur $G$. Dans la terminologie anglo-saxone, on l'appelle ``greatest ambit''. Nous l'appelerons ici tout simplement compactifi\'e de Samuel. Notons que, par le th\'eor\`eme de Gelfand, l'alg\`ebre $C_{\mathbb{C}}(\mathcal{S}(G))$ est isomorphe \`a $E$.
\begin{remark}
Si $G$ est un groupe discret d\'emonbrable, alors le compactifi\'e de Samuel $\mathcal{S}(G)$ co\"\i ncide avec le compactifi\'e de Stone-$\check{C}$ech $\beta G$.
\end{remark}
Soient $G$ un groupe topologique, $Z$ un espace topologique, et $f\colon G\to Z$ une application uniform\'ement continue \`a droite telle que $X:=f(G)$ soit compact; alors $f$ se prolonge en une application $f:\mathcal{S}(G)\longrightarrow X$, encore not\'ee $f$.
D\'efinissons une relation d'\'equivalence sur $\mathcal{S}(G)$ par : $$x\mathcal{R}y\,\,\,\text{si}\,\,\, f(gx) = f(gy)\,\,\,\,\text{pour tout}\,\, g\in G.$$ Il est facile de v\'erifier que le graphe $C$ de $\mathcal{R}$ est ferm\'e dans $\mathcal{S}(G)\times \mathcal{S}(G)$, de sorte que l'espace quotient $X_{f}:=\mathcal{S}(G)/ \mathcal{R}$ est s\'epar\'e. Puisque la surjection canonique $\pi_f:\mathcal{S}(G)\longrightarrow X_{f}$ est continue, l'espace quotient $X_{f}$ est compact. L'application $\bar f:X_{f}\ni [x]\longmapsto f(x)\in X$ est continue, et $f=\bar f\circ \pi_f$.

Soit $F$ une famille d'applications uniform\'ement continues \`a droite de $G$ dans un espace compact $Y$ (vues comme applications $\mathcal{S}(G)\longrightarrow Y$). Posons $X = Y^F$ et soit $f:\mathcal{S}(G)\longrightarrow X$ le produit diagonal de la famille $F$. Notons $X_F = X_f$ dans ce cas-l\`a, et $\pi_F\colon \mathcal{S}(G)\to X_F$ la surjection canonique. Le lemme suivant est imm\'ediat:

\begin{lemma}
Si $F$ s\'epare les points de $\mathcal{S}(G)$ , c'est \`a dire si, pour tous $x,y\in \mathcal{S}(G)$ avec $x\neq y$, il existe $f\in F$ tel que $f(x)\neq f(y)$, alors $\pi_F$ est un hom\'eomorphisme de $\mathcal{S}(G)$ sur $X_F$.
\end{lemma}
\begin{lemma}\label{pif}
Si $F\subseteq F^{\prime}$, alors il existe une surjection continue et $G$-équivariante canonique $\pi^{F^{\prime}}_F$ de $X_{F^{\prime}}$ sur $X_F$.
\end{lemma}
\begin{proof}
Soit $x\in \mathcal{S}(G)$. Pour tout $F\subseteq F^{\prime}$, l'application
$$\pi^{F^{\prime}}_{F}:X_{F^{\prime}}\ni [x]_{\mathcal{R}_{F^{\prime}}}\longmapsto [x]_{\mathcal{R}_{F}}\in X_{F}$$ est bien d\'efinie. De m\^eme, pour tout $F\subseteq F^{\prime}$, on a: $\pi^{F^{\prime}}_{F}\circ \pi_{F^{\prime}}=\pi_{F}$. Donc $\pi^{F^{\prime}}_{F}$ est continue. $G$ op\`ere contin\^ument sur $X_{F}$ par l'action quotient $g[x]_{\mathcal{R}_{F}}=[gx]_{\mathcal{R}_{F}}$ et $\pi^{F^{\prime}}_{F}$  est \'equivariante par rapport \`a cette action.
\end{proof}

\begin{theorem}
Soit $F$ la famille des fonctions à valeurs complexes sur $G$ qui sont uniform\'emenet continues et born\'ees. Alors $X_F$ est l'espace de Gelfand de la $C^{*}$-alg\`ebre engendr\'ee par tous les translat\'es de $F$ dans $RUCB(G) = C_{\mathbb{C}}(\mathcal{S}(G))$.
\end{theorem}
\begin{corollary}
Supposons que $G$ est d\'enombrable, $X$ compact m\'etrisable, et $F$ d\'enombrable. Alors $X_F$ est compact m\'etrisable.
\end{corollary}
Rappelons la notion de syst\`eme projectif de $G$-espaces compacts, en suivant \cite{Bou}.
$G$ \'etant un groupe topologique, un syst\`eme projectif de $G$-espaces est la donn\'ee d'un triplet $(X_{i},\pi_{ij}, I)$ o\`u $(I,\preceq )$ est un ensemble ordonn\'e, $\{X_{i}\}_{i\in I}$ une famille  de $G$-espaces compacts et une famille d'applications \'equivariantes $\pi_{ij}:X_{j}\longrightarrow X_{i}$, pour $i\preceq j$ telle que $\pi_{ii}$ est l'identit\'e de $X_{i}$ pour tout $i\in I$ et $\pi_{i k}=\pi_{i j}\circ \pi_{j k}$ pour tout $i \preceq j \preceq k$. La limite projective $X=\underset{\longleftarrow}{\overset{}{\lim}}\,(X_{i},\pi_{ij})$ du syst\`eme projectif de $G$-espaces $(X_{i},\pi_{ij}, I)$ est le $G$-espace d\'efini comme suit:
si $\underset{i\in I}{\overset{}{\prod}}X_{i}$ est l'espace produit des espaces $(X_{i})_{i\in I}$ et $pr_{i}:\underset{i\in I}{\overset{}{\prod}}X_{i}\longrightarrow X_{i}$ la projection sur le facteur $X_{i}$, alors

 \begin{tabular}{lll}
$\underset{\longleftarrow}{\overset{}{\lim}}(X_{i},\pi_{ij})$
& $=$ &
$\{x=(x_{i})\in \underset{i\in I}{\overset{}{\prod}}X_{i}/\, x_{i}=\pi_{ij}(x_{j})\,\,i\preceq j\}$ \\
& $=$ & $ \{x=(x_{i})\in \underset{i\in I}{\overset{}{\prod}}X_{i}/\, pr_{i}(x)=\pi_{ij}\circ pr_{j}(x)\,\,i\preceq j\}.$ 
\end{tabular}

C'est un sous-espace ferm\'e, donc compact de $\underset{i\in I}{\overset{}{\prod}}X_{i}$. Le groupe $G$ op\`ere contin\^ument sur $\underset{i\in I}{\overset{}{\prod}}X_{i}$ par $g.x=(g.x_{i})$ o\`u $x=(x_{i})$.
La restriction $\pi_{i}$ de la projection $pr_{i}$ \`a X est l'application canonique de $X$ dans $X_{i}$, et $\pi_{i}=\pi_{ij}\circ \pi_{j}$ pour tout $i\preceq j$. On écrit simplement $X=\underset{\longleftarrow}{\overset{}{\lim}}X_{i}$ si aucune confusion n'est possible.

On v\'erifie ais\'ement le lemme:

\begin{lemma}
 Soit $\Phi$ une collection des familles des fonctions de $\mathcal{S}(G)$ dans $X$. Supposons que $\Phi$ est dirig\'ee par inclusion (quels que soient $F,F^{\prime}\in\Phi$, il existe $F^{\prime\prime}\in\Phi$ tel que $F\subseteq F^{\prime\prime}$ et $F^{\prime}\subseteq F^{\prime\prime}$). Alors le triplet $(X_F,\pi^{F^{\prime}}_F,\Phi)$ est un syst\`eme projectif de $G$-espaces.
 \end{lemma}
 \begin{remark}
\label{rem1a}
Supposons que la reunion $\cup \{F\colon F\in\Phi\}$ s\'epare les points de $\mathcal{S}(G)$. Alors la limite projective du syst\`eme $(X_F,\pi^{F^{\prime}}_F)$ est isomorphe, en tant que $G$-espace compact \`a $\mathcal{S}(G)$. En effet, consid\'erons l'application $\Psi:\mathcal{S}(G)\ni x\longmapsto (\pi_F(x))_{F\in \Phi}\in \underset{\longleftarrow}{\overset{}{\lim}}\,X_{F}$. Soit $x,y\in \mathcal{S}(G)$ tel que $x\neq y$. Puisque $\cup \{F\colon F\in\Phi\}$ s\'epare les points de $\mathcal{S}(G)$, il existe $F\in\Phi$ tel que $\pi_F(x)\neq \pi_F(y)$. Donc $\Psi(x)\neq\Psi(y)$ et $\Psi$ est injective. L'application $\Psi$ est surjective par d\'efinition, et \'equivariante.
\end{remark}

\begin{lemma}\label{lem2}
Soient $X$ et $Y$ deux $G$-espaces compacts et soit $\alpha : X\longrightarrow Y$ une application continue et \'equivariante. L'application $\alpha_{\star}: \mathcal{P}(X)\ni \mu\longmapsto \mu\circ \alpha ^{-1} \in \mathcal{P}(Y) $ est \'egalement continue et \'equivariante.
\end{lemma}
\begin{lemma}\label{lemmesurinverse}
S'il existe une mesure de probabilit\'e invariante sur chaque $G$-espace dans un syst\`eme inverse des $G$-espaces compacts, alors il existe une mesure de probabilit\'e invariante sur la limite inverse correspondante.
\end{lemma}
\begin{proof}
Soit $(X_{\alpha},\pi_{\alpha \beta}, I)$ un syst\`eme projectif de $G$-espaces. Notons $X=\underset{\longleftarrow}{\overset{}{\lim}}\,X_{\alpha}$. Par le lemme \ref{lem2}, ce syst\`eme projectif permet d'obtenir une syst\`eme projectif $(\mathcal{P}(X_{\alpha}),(\pi_{\alpha \beta})_{\star}, I)$. Il est clair que $\mathcal{P}(X)=\underset{\longleftarrow}{\overset{}{\lim}}\,\mathcal{P}(X_{\alpha})$. Notons $\mathcal{P}_{inv}(X_\alpha)$ l'espace des mesures de probabilit\'es invariantes sur $X_\alpha$ et consid\'erons le sous-syst\`eme $(\mathcal{P}_{inv}(X_{\alpha}),(\pi_{\alpha \beta})_{\star},I)$ de $(\mathcal{P}(X_{\alpha}),(\pi_{\alpha \beta})_{\star}, I)$. On a \'egalement $\mathcal{P}_{inv}(X)=\underset{\longleftarrow}{\overset{}{\lim}}\,
\mathcal{P}_{inv}(X_{\alpha})$. Puisque $\mathcal{P}_{inv}(X_\alpha)$ est compact et non vide pour tout $\alpha$, on a: $\mathcal{P}_{inv}(X)\neq \emptyset$. Si $\mu\in \mathcal{P}_{inv}(X)$, alors $\mu$ est une mesure de probabilit\'e invariante sur $X$.
\end{proof}

\begin{lemma}
Soit $G$ un sous-groupe dense de $G'$, et $F$ une famille des fonctions de $G$ dans un espace compact $X$. Alors on peut construire $X_F$ \`a partir de $G$, notons-le $X_F(G)$, ainsi qu'\`a partir de $G'$, notons-le $X_F(G')$. Les $G'$-espaces $X_F(G)$ et $X_F(G')$ sont isomorphes entre eux de mani\`ere canonique.
\end{lemma}

\begin{lemma}
Si $X$ est un compact de dimension z\'ero (au sens de Lebesgue), alors $X_F$ est de dimension z\'ero.
\end{lemma}

\begin{remark}
Supposons que $X_F$ contienne un point isol\'e $x_0$; notons
$H=\{g\in G:\,\,gx_0=x_0\}$ le stabilisateur correspondant, qui est un sous-groupe ouvert de $G$. Alors, pour tout $f\in F$ la restriction $f|_{gH}$ est constante: Pour guarantir la non-existence des points isol\'es, nous allons supposer que la famille $F$ v\'erifie la condition $(\star)$ suivante:
  Pour tout voisinage $V$ de $e$, il existe $f\in F$ et $x\in V$ tel que $f(x)\neq f(e)$.
\end{remark}
Dans ces conditions, on a le lemme:
\begin{lemma}
$G$ \'etant un groupe topologique, $X$ un espace compact, et $X_F$ d\'efini comme pr\'ec\'edement avec $F$ v\'erifiant la condition $(\star)$, alors ou bien $X_F$ est fini, ou bien $X_F$ ne contient aucun point isol\'e.
\end{lemma}

\begin{corollary}\label{corcantor}
Si $G$ est infini d\'enombrable, $X$ compact de dimension z\'ero et $F$ d\'enombrable et v\'erifie $(\star)$, alors $X_F$ est hom\'eomorphe \`a l'espace de Cantor.
\end{corollary}

\begin{definition}
Un groupe topologique $G$ est dit {\em non archim\'edien} s'il est s\'epar\'e et poss\`ede une base de voisinages de l'\'el\'ement neutre form\'e des sous-groupes ouverts.
\end{definition}

L'ensemble des groupes polonais non archimédiens comprend:
\begin{enumerate}
  \item le groupe sym\'etrique infini $S_{\infty}$ de toutes les bijections de $\mathbb{N}$ dans $\mathbb{N}$, muni de la topologie de la convergence simple,
  \item le groupe $Homeo(D^{\aleph_{0}})$ des hom\'eomorphismes de l'ensemble de Cantor muni de la topologie de la convergence uniforme,
  \item les groupes localement compacts totalement discontinus (voir \cite{Bou}, chapitre III, section $4$, no. $6$).
\end{enumerate}

Les groupes polonais non archim\'ediens jouent un r\^ole important en logique o\`u ils sont les groupes des automorphismes des structures de Fraïss\'e \cite{becker}.

\begin{theorem}(\cite{vp2})
Un groupe topologique $G$ est non archim\'edien si et seulement si les applications continues de $\mathcal{S}(G)$ dans $D = \{0,1\}$ s\'eparent les points de $\mathcal{S}(G)$.
\end{theorem}
On d\'eduit de tout ce qui pr\'ec\`ede le th\'eor\`eme fondamental suivant:
\begin{theorem}\label{coroinverse}
Si $G$ est non archim\'edien et polonais, alors $\mathcal{S}(G)$ se d\'eveloppe en limite projective d'un syst\`eme des $G$-espaces compacts m\'etrisables de dimension z\'ero. Si $G$ est de plus infini, alors $\mathcal{S}(G)$ est la limite d'un tel système dont tous les $G$-espaces sont homéomorphes \`a l'espace de Cantor.
\end{theorem}

\section{Espaces test pour la moyennabilit\'e}

\begin{definition}
Pour une partie compacte $X\subset E$ d'un espace localement convexe
$E$, on appelle {\em barycentre} d'une mesure $\mu\in \mathcal{P}(X)$ une forme lin\'eaire $b(\mu)\in E^{\prime\prime}$, d\'efinie par l'\'egalit\'e
$b(\mu)(\varphi)=\mu(\varphi \mid X)$.
\end{definition}

\begin{remark}
Si $X$ est un compact convexe de $E$, alors il existe un unique barycentre $b(\mu)$ pour tout $\mu\in \mathcal{P}(X)$, et de plus $b(\mu)\in X$ (c.à.d. $b(\mu)$ est une mesure de Dirac); voir \cite{bour} Chapitre 4, $\S$ 7). Dans ce cas, l'application $b=b_{X}: \mathcal{P}(X) \longrightarrow X$ est bien d\'efinie, continue (voir \cite{bour} Chapitre 3, $\S$ 3)  et \'evidemment affine.
\end{remark}

\begin{lemma} \cite{boga}\label{equi}
 Si $X\subset E$ est convexe et le groupe $G$ op\`ere affinement sur
$X$, alors l'application $b_{X}: \mathcal{P}(X) \longrightarrow X$ est \'equivariante par rapport \`a l'action de $G$ sur $\mathcal{P}(X)$.
\end{lemma}

Par cons\'equent, on a le r\'esultat suivant \cite{boga}:

\begin{theorem}\label{th1}
Pour une action affine d'un groupe topologique $G$ sur un espace compact convexe, les propositions suivantes sont \'equivalentes:
\begin{enumerate}
  \item l'action poss\`ede un point fixe,
  \item l'action poss\`ede une mesure de probabilit\'e bor\'elienne invariante.
\end{enumerate}
\end{theorem}

\begin{lemma}\label{pr2}
Un groupe topologique $G$ est moyennable si et seulement s'il existe une mesure de probabilit\'e invariante sur son compactifi\'e de Samuel $\mathcal{S}(G)$.
\end{lemma}

\begin{lemma}(\cite{kpt})\label{poloinverse}
Si $G$ est un groupe polonais, alors il existe un syst\`eme projectif de $G$-espaces compacts et m\'etrisables
$(X_{\alpha},\pi_{\alpha \beta}, I)$ tel que $\mathcal{S}(G)=\underset{\longleftarrow}{\overset{}{\lim}}\,X_{\alpha}$.
\end{lemma}

Le r\'esultat suivant est bien connu au moins pour la classe des groupes discrets d\'enombrables. Nous rappelons l'argument quand m\^eme.

\begin{proposition}\label{pr4}
Un groupe polonais $G$ est moyennable si et seulement si toute action continue de $G$ sur un espace compact et m\'etrisable, poss\`ede une mesure de probabilit\'e invariante.
\end{proposition}

\begin{proof}
La n\'ecessit\'e est \'evidente. Montrons la suffisance. Par le lemme \ref{poloinverse}, $\mathcal{S}(G)=\underset{\longleftarrow}{\overset{}{\lim}}\,X_{\alpha}$. $X_\alpha$ \'etant un $G$-espace compact et m\'etrisable pour tout $\alpha$. Par hypoth\`ese, il existe sur chaque $G$-espace compact m\'etrisable $X_\alpha$ une mesure de probabilit\'e invariante. Donc il existe une mesure de probabilit\'e invariante sur $S(G)$ par le lemme \ref{lemmesurinverse} et $G$ est moyennable.
\end{proof}

\begin{remark}\label{rem1}
\'Evidement, dans la démonstration de la proposition \ref{pr4}, on peut supposer sans perte de g\'en\'eralit\'e que tous les $G$-espaces $X$ sont infinis.
\end{remark}

\begin{theorem}\label{main1}
Un groupe polonais $G$ est moyennable si et seulement si toute action continue de $G$ sur le cube de Hilbert $I^{\aleph_{0}}$ poss\`ede une mesure de probabilit\'e bor\'elienne invariante.
\end{theorem}

\begin{proof}La n\'ecessit\'e est \'evidente. Montrons la suffisance.
Soit $X$ un $G$-espace compact et m\'etrisable. Par la remarque \ref{rem1}, on peut supposer que $X$ est infini. $\mathcal{P}(X)$ est donc un sous-espace compact m\'etrisable de dimension infinie de $\mathbb{R}^{\mathcal{C}(X)}$. Par le th\'eor\`eme de Keller (voir \cite{bes}), $\mathcal{P}(X)$ est hom\'eomorphe au cube de Hilbert $I^{\aleph_{0}}$. Ainsi l'action de $G$ sur $\mathcal{P}(X)$ poss\`ede une mesure de probabilit\'e bor\'elienne invariante. L'action de $G$ sur $\mathcal{P}(X)$ \'etant affine, elle poss\`ede par le th\'eor\`eme \ref{th1} un point fixe $\mu \in \mathcal{P}(X)$, qui est une mesure de probabilit\'e bor\'elienne invariante pour l'action initiale de $G$ sur $X$.
\end{proof}

\begin{remark}
L'id\'ee d'utiliser le th\'eor\`eme de Keller dans le contexte dynamique appartient de toute \'evidence \`a Uspenskij, qui \'etait  le premier \`a l'employer dans \cite{usp}.
\end{remark}

\begin{theorem}
 L'ensemble de Cantor $D^{\aleph_{0}}$ est un espace test pour la moyennabilit\'e des groupes polonais non archimédiens. Autrement dit, un groupe polonais non archimédien $G$ est moyennable si et seulement si toute action continue de $G$ sur $D^{\aleph_{0}}$ poss\`ede une mesure de probabilit\'e invariante.
\end{theorem}

\begin{proof}
La n\'ecessit\'e est \'evidente. Montrons la suffisance. $G$ \'etant un groupe polonais non archimédien, il existe par le th\'eor\`eme \ref{coroinverse}, un syst\`eme projectif de $G$-espaces $(X_{\alpha},\pi_{\alpha \beta}, I)$ avec $X_{\alpha}\cong D^{\aleph_{0}}$ pour tout $\alpha \in I$ tel que $\mathcal{S}(G)=\underset{\longleftarrow}{\overset{}{\lim}}\,X_{\alpha}$. Par hypoth\`ese, il existe sur chaque $G$-espace $X_{\alpha}$ une mesure de probabilit\'e invariante $\mu_\alpha$. Par le lemme \ref{lemmesurinverse}, il existe une mesure de probabilit\'e invariante $\mu$ sur $\mathcal{S}(G)$ et $G$ est moyennable.
\end{proof}

\begin{remark}
Pour la th\'eorie des groupes moyennables localement compacts, voir \cite{jp}, \cite{paterson}, \cite{boga}. Pour celle des groupes moyennables non localement compact, voir \cite{dlH}.
\end{remark}

\section{Sur les espaces test pour la moyennabilit\'e extr\^eme}

\begin{definition}
Un groupe topologique $G$ est dit {\em extr\^emement moyennable} si toute action continue de $G$ sur un espace compact $K$ poss\`ede un point fixe.
\end{definition}
Des exemples des groupes extr\^emement moyennables sont nombreux et ils comprennent:
\begin{enumerate}
   \item Le groupe unitaire $\mathcal{U}(\ell^{2})$, muni de la topologie forte (Gromov et Milman \cite{gromov}).
   \item Le groupe $Aut\,(X,\mu)$ des automorphismes mesurables pr\'eservant la mesure $\mu$ d'un espace borelien $(X,\mu)$ muni de la topologie faible (Giordano et Pestov \cite{giopes1}).
   \item Le groupe $Aut\,(\mathbb{Q},\leq)$ des bijections de $\mathbb{Q}$ dans lui-m\^eme qui pr\'eservent l'ordre muni de la topologie de la convergence simple (Pestov \cite{vp2}).
 \end{enumerate}

Le r\'esulat suivant est imm\'ediat.

\begin{proposition}\label{pr3}
Un groupe topologique $G$ est extr\^emement moyennable si et seulement si l'action de $G$ sur son compactifi\'e de Samuel $\mathcal{S}(G)$ poss\`ede un point fixe.
\end{proposition}

\begin{theorem}\label{main3}
Un groupe polonais non archim\'edien $G$ est extr\^emement moyennable si et seulement si toute action continue de
$G$ sur l'ensemble de Cantor $D^{\aleph_{0}}$ poss\`ede un point fixe.
\end{theorem}

\begin{proof}
  La n\'ecessit\'e est \'evidente. Montrons la suffisance. Montrons que l'action canonique de $G$ sur $\mathcal{S}(G)$ poss\`ede un point fixe. Comme
  $G$ est polonais et non archimedien, il existe par le th\'eor\`eme \ref{coroinverse}, un syst\`eme projectif de $G$-espaces $(X_{\alpha},\pi_{\alpha \beta}, I)$ avec $X_{\alpha}\cong D^{\aleph_{0}}$ pour tout $\alpha \in I$ tel que $\mathcal{S}(G)=\underset{\longleftarrow}{\overset{}{\lim}}\,X_{\alpha}$. Par hypoth\`ese, il existe sur chaque $G$-espace $X_{\alpha}$ un point fixe $x_{\alpha}$. Notons $\pi_{\alpha}$ la restriction de la projection $pr_{\alpha}$ \`a $\mathcal{S}(G)=\underset{\longleftarrow}{\overset{}{\lim}}\,X_{\alpha}$ et posons $M_{\alpha}=\pi_{\alpha}^{-1}(x_{\alpha})$. Les applications $\pi_{\alpha}$ \'etant surjectives, on a $M_{\alpha}\neq \emptyset$ pour tout $\alpha$. La famille $(M_{\alpha})_{\alpha \in I}$ est centr\'ee car $x=(x_{\alpha_{1}},...,x_{\alpha_{n}})\in \underset{i=1}{\overset{n}{\cap}}\,M_{\alpha_{i}}$ pour tout $i=1,2,..,n$. L'espace
   $\mathcal{S}(G)$ \'etant compact, $\underset{\alpha \in I}{\overset{}{\bigcap}}\,M_{\alpha}\neq\emptyset$. Tout point de $\underset{\alpha \in I}{\overset{}{\bigcap}}\,M_{\alpha}$ est fixe pour l'action continue de $G$ sur $\mathcal{S}(G)$.
\end{proof}
\begin{question}

Existe-t-il un espace test pour les groupes polonais extr\^emement moyennables?
\end{question}

Le th\'eor\`eme du point fixe de Schauder affirme que toute fonction continue de $I^{\aleph_{0}}$ dans $I^{\aleph_{0}}$ poss\`ede un point fixe. En particulier, toute action continue du groupe discret $\mathbb{Z}$ sur $I^{\aleph_{0}}$ par hom\'eomophismes poss\`ede un point fixe. Ceci permet de conclure que le cube de Hilbert $I^{\aleph_{0}}$ ne peut pas \^etre un espace test pour la moyennabilit\'e des groupes polonais. En effet, le th\'eor\`eme de Ellis \cite{ellis} affirme que tout groupe discret agit librement sur un espace compact et par cons\'equent, n'est pas extr\^emement moyennable.\\
On peut n\'eanmoins observer qu'il existe un espace test compact s\'eparable non n\'ecessairement m\'etrisable pour les groupes polonais extr\^emement moyennables.
En effet, notons par $\mathcal{P}_{0}$ l'ensemble de tous les groupes polonais non-extr\^emement moyennables deux \`a deux non-isomorphes et choisissons pour tout $G\in \mathcal{P}_{0}$ un $G$-espace compact et m\'etrisable $X_{G}$ sans points fixes. L'espace $X=\underset{G\in \mathcal{P}_{0}}{\overset{}{\prod}}X_{G}$ est un espace test s\'eparable compact (non n\'eccessairement m\'etrisable) pour la moyennabilit\'e extr\^eme des groupes polonais .
Il est clair que l'action produit de $G$ sur $X$ est continue et sans points. Pour conclure, et grâce au c\'el\`ebre th\'eor\`eme de Hewitt \cite{hewitt} et Pondiczery \cite{pondi}, il suffit de montrer que
 $|\mathcal{P}_{0}|\leq 2^{\aleph_0}$. Notons $F_{\infty}$ le groupe libre avec un nombre infini d\'enombrable de g\'en\'erateurs. Notons $\mathcal{P}$ l'ensemble des groupes  polonais et $\mathcal{D}$ l'ensemble de toutes les pseudo-m\'etriques sur $F_{\infty}$. Il est clair que $|\mathcal{D}|\leq |\mathbb{R}^{\mathbb{Z}}|=2^{\aleph_0}$.
Montrons que $|\mathcal{P}| \leq |\mathcal{D}|$.\\
Soit $d$ une pseudo-m\'etrique sur $F_{\infty}$ invariante \`a gauche. $H_{d}=\{x\in F_{\infty}, d(x,e)=0\}$ est un sous-groupe de $F_{\infty}$. La distance d\'efinie sur $F_{\infty}/H_{d}$ par $\widehat{d}(xH_{d},yH_{d})=d(x,y)$ est invariante par translation \`a gauche.
Notons $G_{d}$ le complet\'e de l'espace m\'etrique $(F_{\infty}/H_{d},\widehat{d})$. Si  $H_{d}$ est un sous-groupe normale de $F_{\infty}$, le groupe topologique $G_{d}$ est un groupe polonais et chaque groupe polonais est de la forme $G_{d}$. Notons par $\mathcal{D}_{N}$ le sous-  ensemble de $\mathcal{D}$ constitu\'e des pseudo-m\'etriques $d$ telles que $H_{d}$ soit normale. Nous avons donc une application surjective: $\mathcal{D}_{N}\ni d\longmapsto G_{d}\in \mathcal{P}$. Ainsi $|\mathcal{P}|\leq |\mathcal{D}_{N}|\leq |\mathcal{D}|$.\\

Rappelons qu'un groupe topologique $G$ est dit {\em monoth\'etique}, s'il existe un sous-groupe $H$ de $G$ qui est \`a la fois cyclique et dense, et {\em sol\'enoïde}, s'il existe un homomorphisme continu $f$ de $\mathbb{R}$ dans $G$ dont l'image est partout dense dans $G$.

\begin{remark}
\begin{enumerate}
  \item Il est clair que tout groupe monoth\'etique ou sol\'enoide est ab\'elien, donc moyennable.
  \item Il est facile de d\'eduire du th\'eor\`eme du point fixe de Schauder que toute action continue d'un groupe monoth\'etique ou d'un groupe sol\'enoïde sur le cube de Hilbert $I^{\aleph_{0}}$ admet un point fixe.
\end{enumerate}
\end{remark}
\begin{question}
 Toute action continue d'un groupe polonais moyennable sur le cube de Hilbert $I^{\aleph_{0}}$ poss\`ede-t-elle un point fixe?
   M\^eme question pour un groupe moyennable discret.
\end{question}

\section{Espaces test pour la moyennabilit\'e \`a l'infini}

Gr\^ace \`a la condition de Reiter ($P_{1}$) \cite{paterson}, un groupe discret d\'emonbrable $G$ est moyennable si et seulement si l'action triviale de $G$ sur $\{\star\}$ est moyennable.
La moyennabilit\'e \`a l'infini est plus g\'en\'erale que la moyennabilit\'e. En effet, le groupe libre \`a deux g\'en\'erateurs $F_{2}$ est moyennable \`a l'infini (\cite{BO}, prop.$5.1.8$) mais n'est pas moyennable (\cite{paterson}, p. $6$).

Le r\'esultat suivant est bien connu de la th\'eorie de la moyennabilit\'e \`a l'infini.

\begin{theorem}(Assertion $1$ dans \cite{dran})
Un groupe d\'enombrable $G$ admet une action moyennable sur un espace compact et m\'etrisable si et seulement si son action sur son compactifi\'e de Stone-$\check{C}$ech  $\beta G$ est moyennable.
\end{theorem}

\begin{lemma}\label{lemmoytopo}
Soit $G$ un groupe d\'enombrable moyennable \`a l'infini. Notons $(b_n)$ la suite des applications correspondantes de $\mathcal{S}(G)$ dans $\mathcal{P}(G)$ et posons $F= \{b_n: n\in \mathbb{N}\}$. Alors l'action de $G$ sur $X_F$ est moyennable \`a l'infini.
\end{lemma}

\begin{proof}
$G$ \'etant moyennable \`a l'infini, il existe une suite d'applications continues $b^{n}:\mathcal{S}(G)=\beta G\longrightarrow \mathcal{P}(G)$ telles que: $\underset{n\longrightarrow \infty}{\overset{}{\lim}}\,\,\underset{x\in \beta G}{\overset{}{\sup}}\|gb^{n}_{x}-b^{n}_{gx}\|_{1}=0$ pour tout $g\in G$.
Pour tout $g\in G$, notons $\overline{g}:\beta G \ni x\longmapsto gx \in \beta G$ l'hom\'eomorphisme de $\beta G$ sur lui m\^eme produit par $g$. Consid\'erons le produit diagonal
\[f=\Delta_{(g,n)\in G\times \mathbb{N}}(b_{n}\circ \overline{g}):\beta G\longrightarrow (\mathcal{P}(G))^{G\times \mathbb{N}}\]
d\'efini par $f(x)=(b_{gx}^{n})_{(g,n)\in G\times \mathbb{N}}$.
Il est clair que $f$ est continue. La relation d'\'equivalence $\mathcal{R}_{F}$ est d\'efinie sur $\beta G$ par: $(x,y)\in \mathcal{R}_{F}\Longleftrightarrow b_{gx}^{n}=b^{n}_{gy}$ pour tout $n\in \mathbb{N}$ et $g\in G$.
   Notons encore par $f$ l'application $f:\beta G\longrightarrow f(\beta G)$, et par $\pi_F:\beta G\longrightarrow \beta G/ \mathcal{R}_{F}$ la surjection canonique. Il existe une application continue $\bar f$ tel que: $f=\bar f\circ \pi_F$. Posons $X_F=\beta G/ \mathcal{R}_{F}$.
    Consid\'erons l'application $\widetilde{b}^{n}:f(\beta G)\longrightarrow \mathcal{P}(G)$ d\'efinie par $\widetilde{b}^{n}=\pi_{e,n}$ o\`u $\pi_{e,n}$ est la restriction de la projection \`a $f(\beta G)$. Ainsi, les applications $c^{n}:\widetilde{b}^{n}\circ \bar f:X_F\longrightarrow \mathcal{P}(G)$ sont continues.

Soit $g\in G$, on a:

    \begin{tabular}{lll}
$\underset{[x]\in X_F }{\overset{}{\sup}}\|gc^{n}_{[x]}-c^{n}_{g[x]}\|_{1}$
& $=$ &
$\underset{[x]\in X_F }{\overset{}{\sup}}\|g(\widetilde{b}^{n}\circ \bar f)_{[x]}-(\widetilde{b}^{n}\circ \bar f)_{g[x]}\|_{1}$ \\
& $=$ & $ \underset{[x]\in X_F }{\overset{}{\sup}}\|g(\widetilde{b}^{n}(\bar f([x])))-\widetilde{b}^{n}(\bar f(g[x]))\|_{1}$ \\
& $= $ & $ \underset{x\in\beta G}{\overset{}{\sup}}\|g(\widetilde{b}^{n}(\bar f([x])))-\widetilde{b}^{n}(\bar f([gx]))\|_{1}$ \\
& $=$ & $\underset{x\in\beta G}{\overset{}{\sup}}\|g(\widetilde{b}^{n}(f(x)))-\widetilde{b}^{n}(f(gx))\|_{1}$\\
& $=$ & $\underset{x\in\beta G}{\overset{}{\sup}}\|g(\widetilde{b}^{n}(b_{hx}^{n}))_{(h,n)\in G\times \mathbb{N}})-\widetilde{b}^{n}((b_{hgx}^{n}))_{(h,n)\in G\times \mathbb{N}}\|_{1}$\\
& $=$ & $\underset{x\in\beta G}{\overset{}{\sup}}\|gb_{x}^{n}-b^{n}_{gx}\|_{1}$
\end{tabular}

Ainsi, $\underset{n\longrightarrow \infty}{\overset{}{\lim}}\,\underset{[x]\in X_{F} }{\overset{}{\sup}}\|gc^{n}_{[x]}-c^{n}_{g[x]}\|_{1}=0$.
\end{proof}

\begin{theorem} Un groupe discret d\'emonbrable $G$ est moyennable \`a l'infini si et seulement s'il admet une action moyennable sur l'ensemble de Cantor $D^{\aleph_{0}}$.
\end{theorem}

\begin{proof}
La suffisance est \'evidente. Montrons la necessit\'e. Si $G$ poss\`ede une action moyennable sur son compactifi\'e de Stone-$\check{C}$ech, alors il existe une suite d'applications $b^{n}:\beta G\longrightarrow \mathcal{P}(G)$ comme dans la d\'efinition de la moyennabilit\'e \`a l'infini.
Pour tout $g\in G$, notons $\overline{g}:\beta G\ni x\longmapsto gx \in \beta G$ l'hom\'eomorphisme de $\beta G$ sur lui m\^eme. L'espace $b^{n}(\beta G)$ \'etant compact et m\'etrisable, il existe une surjection continue $f^{n}:D^{\aleph_{0}}\longrightarrow b^{n}(\beta G)$.  Soit $g\in G$, alors $b_{g}^{n}\in b^{n}(\beta G)$. Puisque $f^{n}$ est surjective, il existe $c_{g}\in D^{\aleph_{0}}$ tel que $f^{n}(c_{g})=b_{g}^{n}$. On a ainsi une application $T^{n}:G\ni g\longmapsto c_{g}\in D^{\aleph_{0}}$. Cette application se prolonge de mani\`ere unique en une application continue $\beta T^{n}:\beta G\longmapsto D^{\aleph_{0}}$.
Pour tout $g\in G$, on a:
$(f^{n}\circ \beta T^{n})(g)=f^{n}(\beta T^{n}(g))=f^{n}(c_{g})=b_{g}^{n}$.
Ainsi, $f^{n}\circ \beta T^{n}=b^{n}$ sur $G$. Puisque $G$ est dense dans $\beta G,\,\,\,f^{n}\circ \beta T^{n}=b^{n}$ sur $\beta G$.
Posons: $c^{n}=\beta T^{n}\circ \overline{g}:\beta G\longrightarrow \{0,1\}^{\aleph_{0}},\,\,F= \{c^{n}: n\in \mathbb{N}\}$ et $X=D^{\aleph_{0}}$. Alors $X_{F}\cong D^{\aleph_{0}}$ et l'action de $G$ sur $X_F$ est moyennable \`a l'infini par le lemme \ref{lemmoytopo}.
\end{proof}

\begin{theorem}\label{corohilbert}
Un groupe discret d\'emonbrable $G$ est moyennable \`a l'infini si et seulement s'il admet une action moyennable sur le cube de Hilbert $I^{\aleph_{0}}.$
\end{theorem}

Faisons d'abord les rappels suivants:

\begin{lemma}(Lemme $3.6$ dans \cite{higson})\label{lemhig}
Un groupe discret d\'enombrable $G$ est moyennable \`a l'infini si et seulement s'il existe une suite d'applications $b^{n}:G\longrightarrow \mathcal{P}(G)$ tel que:
\begin{enumerate}
  \item pour tout $n$, il existe $F_{n}\subset G$ fini tel que $supp(b^{n}_{g})\subset F_{n}$ pour tout $g\in G$;
  \item $\underset{n\longrightarrow \infty}{\overset{}{\lim}}\,\,\underset{h\in G}{\overset{}{\sup}}\|gb^{n}_{h}-b^{n}_{gh}\|_{1}=0$ pour tout $g\in G.$
\end{enumerate}
\end{lemma}

\begin{lemma}\label{lemcantor}
Soit $G$ un groupe discret d\'emonbrable. Si l'action de $G$ sur l'ensemble de Cantor $D^{\aleph_{0}}$ est moyennable \`a l'infini, alors les images des $b^{n}$ sont finies.
\end{lemma}

\begin{proof}
Soit $b^{n}:D^{\aleph_{0}}\longrightarrow \mathcal{P}(G)$ une suite d'applications comme dans la d\'efinition de la moyennabilit\'e \`a l'infini.
Puisque la dimension de Lebesgue de l'espace de Cantor $D^{\aleph_{0}}$ est z\'ero, il existe une partition finie $\gamma=(A_{i})_{i=1,2,..,n}$ de $D^{\aleph_{0}}$ en sous-ensembles ouverts et ferm\'es, telle que l'image par $b^{n}$ de chaque \'el\'ement de $\gamma$ est contenue dans une des boules $B_{\frac{1}{n}}$ du rayon $1/n$ dans $\ell^1(G)$. Notons $c_{i}$ les centres des boules correspondantes.
  Pour tout $A_{i}\in\gamma$, consid\'erons l'application $(b^{n})^\prime$ d\'efinie par :
   $(b^{n})_{x}^\prime=c_{i}\,\, \text{si}\,\, x\in A_{i}\in \gamma$. Soit $x\in D^{\aleph_{0}}$, il existe par d\'efinition de $b^{n}$ et $(b^{n})^{\prime}$, une boule $B_{\frac{1}{n}}$ telle que:
\[\|b^n_x-(b^n)'_x\|_{1}\leq \underset{p,q \in B_{\frac{1}{n}}}{\overset{}{\sup}}\|p-q\|_{1}=diam(B_{\frac{1}{n}})=\frac{2}{n}.\]
         D'o\`u $\underset{n\longrightarrow \infty}{\overset{}{\lim}}\,\,\underset{x\in D^{\aleph_{0}} }{\overset{}{\sup}}\|g(b^{n})^{\prime}_{x}-(b^{n})^{\prime}_{gx}\|_{1}=0.$
\end{proof}

\begin{proof}[D\'emonstration du Th\'eor\`eme \ref{corohilbert}]
La suffisance est \'evidente. Montrons la necessit\'e.
Si $G$ est un groupe discret et d\'enombrable admettant une action moyennable sur l'ensemble de Cantor, alors il existe une suite d'applications $b^{n}:D^{\aleph_{0}}\longrightarrow \mathcal{P}(G)$ comme dans la d\'efinition de la moyennabilit\'e \`a l'infini.
 Par le lemme \ref{lemhig}, on peut supposer sans perte de g\'en\'eralit\'e que pour tout $n$, il existe $F_n \subseteq G$ fini tel que $supp(b^{n}_{x})\subset F_{n}$ pour tout $x\in D^{\aleph_{0}}$. Par le lemme \ref{lemcantor}, supposons que les images des $b^n$ sont finies et notons $(c_{i})_{i\in I=1,2,..,n}$) les images de toutes les applications $b^n$.
  Posons $A_{i}=(b^{n})_{c_{i}}^{-1}$ et consid\'erons
\[c^{n}:\mathcal{P}(D^{\aleph_{0}})\ni \mu\longmapsto \underset{i=1}{\overset{n}{\Sigma}} \mu(A_{i}) c_{i}\in \mathcal{P}(G).\]
     Les applications $c^{n}$ sont clairement affines et continues par rapport \`a la topologie vague sur $P(D^{\aleph_{0}})$. Ce sont pr\'ecisement les prolongements affines des applications $b^{n}$ sur $\mathcal{P}(D^{\aleph_{0}})$.
     Notons $\mathcal{P}_{0}$ le sous-espace de $\mathcal{P}(D^{\aleph_{0}})$ form\'e des mesures \`a support fini.
     Soit $\mu=\underset{i=1}{\overset{n}{\Sigma}}\alpha_{i}\delta_{x_{i}}\in \mathcal{P}_{0}$, on a: $c^{n}_{\mu}=\underset{i=1}{\overset{n}{\Sigma}}\alpha_{i}c^{n}_{\delta_{x_{i}}}$.
     Or $c^{n}_{\delta_{x_{i}}}=\underset{j=1}{\overset{n}{\Sigma}}\delta_{x_{i}}(A_{j})c_{j}=c_{i}=b^{n}_{x_{i}}$. D'o\`u $c^{n}_{\mu}=\underset{i=1}{\overset{n}{\Sigma}}\alpha_{i}b^{n}_{x_{i}}$. De m\^eme, $c^n_{g\mu}=\Sigma_i \alpha_i b^n_{g x_i}$ car $g\mu=g\underset{i=1}{\overset{n}{\Sigma}}\alpha_{i}\delta_{x_{i}}=\underset{i=1}{\overset{n}{\Sigma}}\alpha_{i}\delta_{gx_{i}}$.\\
      Pour tout $ \varepsilon> 0$,
\[\|gc^n_\mu - c^n_{g\mu} \|_{1}= \|\Sigma_i \alpha_i gb^n_{x_i} - \Sigma_i \alpha_i b^n_{g x_i}\|_{1}\leq \Sigma_i \alpha_i \|gb^n_{x_i} - b^n_{g x_i}\|_{1} \leq \Sigma_i \alpha_i \varepsilon = \varepsilon.\]
      $\mathcal{P}_{0}$ \'etant dense dans $\mathcal{P}(D^{\aleph_{0}})$ et les applications $c^{n}$ continues, on a:
$\underset{\mu \in \mathcal{P}(D^{\aleph_{0}})}{\overset{}{\sup}}\|gc^{n}_{\mu}-c^{n}_{g\mu}\|_{1}<\varepsilon$. D'o\`u $\underset{n\longrightarrow \infty}{\overset{}{\lim}}\,\underset{\mu \in \mathcal{P}(D^{\aleph_{0}})}{\overset{}{\sup}}\|gc^{n}_{\mu}-c^{n}_{g\mu}\|_{1}=0$.

On conlut que l'action de $G$ sur $\mathcal{P}(D^{\aleph_{0}})$ est moyennable. Par le th\'eor\`eme de Keller, $\mathcal{P}(D^{\aleph_{0}})$ est hom\'eomorphe \`a $I^{\aleph_{0}}$.
\end{proof}
Remerciements. Les auteurs sont reconnaissants à Pierre de la Harpe pour ses nombreuses remarques utiles sur une version originale de cette note. Nous remercions le CRSNG Canada (V.G.P., B.R.M.), le MAÉCI canada (B.R.M), et le gouvernement de Libye (Y.A.-G.) pour leur support de ce projet de recherche.

\end{document}